\documentclass[multi]{cambridge7A}
\usepackage[UKenglish]{babel}
\usepackage[longnamesfirst,sectionbib]{natbib}
\usepackage{chapterbib}
\usepackage{amsmath,amssymb,amsthm,amsfonts}
\usepackage{enumerate,url}
\usepackage{graphicx}

\setcitestyle{authoryear,round,semicolon}

\theoremstyle{plain}
\newtheorem{theorem}{Theorem}[section]

\allowdisplaybreaks[1]

\newcommand{\mbR}{\mathbb{R}}
\newcommand{\tP}{\mathbf {P}}
\newcommand{\tE}{\mathbf {E}}
\newcommand{\cS}{\cal S}

\newcommand{\Fd}[1]{{#1}^{\ast}}
\newcommand{\swFd}[1]{{#1}^{\ast \!\!\!\circ}}
\newcommand{\sweep}[1]{{#1}^{{\circ}}}
\newcommand{\rBo}[1]{{B}^{(#1)}}
\newcommand{\rB}[2]{{B}^{(#2)}_{#1}}

\newcommand{\cv}[2]{{\mathfrak C}[#1,#2]}

\providecommand*\Index[1]{#1\index{#1}}
\providecommand*\undex[1]{} 

\begin{document}
\alphafootnotes
\author[J. D. Biggins]{J. D. Biggins\footnotemark }
\chapter{Branching out}
\footnotetext{Department of Probability \& Statistics, Hicks Building,
  University of Sheffield, Sheffield S3 7RH; J.Biggins@sheffield.ac.uk}
\arabicfootnotes
\contributor{John D. Biggins
  \affiliation{University of Sheffield}}
\renewcommand\thesection{\arabic{section}}
\numberwithin{equation}{section}
\renewcommand\theequation{\thesection.\arabic{equation}}
\numberwithin{figure}{section}
\renewcommand\thefigure{\thesection.\arabic{figure}}

\begin{abstract}
Results on the behaviour of the rightmost particle in the $n$th generation in the branching random walk are reviewed and 
the phenomenon of anomalous spreading speeds, noticed recently
in related deterministic models, is considered.  The relationship between such results and certain coupled reaction-diffusion equations is indicated.
\end{abstract}

\subparagraph{AMS subject classification (MSC2010)}60J80

\section{Introduction}

I arrived at the University of Oxford in the autumn of 
1973 for postgraduate study. My intention at
that point was to work in Statistics\index{statistics}. The first year of study was a mixture of 
taught courses and designated reading on three areas (Statistics, Probability, and Functional Analysis, in my case) in the ratio 2:1:1 and a dissertation on the main area.  
As part of the Probability\index{probability|(} component, I attended a
graduate course that was an exposition, by its author\index{Hammersley,
J. M.|(}, of the material in \citet{MR0370721}, which had grown out of his
contribution to the discussion of John's\index{Kingman, J. F. C.!influence|(} invited paper on subadditive\index{subadditivity} ergodic theory \citep{MR0356192}.
A key point of Hammersley's contribution  was that the postulates used did not
cover the time to the first birth in the $n$th generation in a Bellman--Harris
process\index{Bellman, R. E.!Bellman--Harris
process}.\footnote{Subsequently, \citet{MR806224} established the theorem
under weaker postulates.}\index{Liggett, T. M.|pagenote}  \citet{MR0370721} showed, among other things, that these quantities did indeed exhibit the anticipated limit behaviour in probability.
I decided not to be examined on this course, which was I believe a wise decision,  but I was intrigued by the material. 
That interest 
turned out to be critical a few months later. 

By the end of the academic year I had concluded that
I wanted to pursue research in Probability\index{probability|)} rather than Statistics and 
asked to have John as supervisor. He agreed. Some time later we met and
he asked me whether I had any particular interests already---I
mentioned Hammersley's lectures. 
When I met him he was in the 
middle of preparing something
(which I could see, but not read upside down). 
He had what seemed to be a pile of written pages, a part written page and a pile of blank paper.
There was nothing else on the desk. 
A few days later a photocopy of 
a handwritten version of  \citet{MR0400438}\index{branching random walk (BRW)|(}\index{branching process!age-dependent branching process|(}, essentially identical to the published version,  appeared in my pigeon-hole with
the annotation ``the multi\-type version is an obvious
problem''\index{branching process!multitype age-dependent branching process}---I am sure this document was what he was writing when I saw him. 
(Like all reminiscences, this what I recall, but it is not necessarily what happened.) This set me going. 
For the next two years, it was a privilege to have John as my thesis supervisor. He supplied exactly what I needed at the time: an initial sense of direction, a strong encouragement to independence, an occasional nudge on the tiller about what did or did not seem tractable, the discipline of explaining orally what I had done, and a ready source on what was known, and where to look for it. However, though important, none of these get to the heart of the matter, which is that I am particularly grateful to have had that period of contact with, and opportunity to appreciate first-hand, such a  gifted mathematician\index{Kingman, J. F. C.!influence|)}.

\citet{MR0400438} considered the problem Hammersley\index{Hammersley, J. M.|)} had raised in its own right, rather than
as an example of, and adjunct to, the general theory of subadditive processes.
Here, I will say something about some recent significant 
developments on the first-birth problem\undex{branching random walk (BRW)!first-birth problem}. I will also go back to my beginnings, 
by outlining something new about the multi\-type version 
that concerns the phenomenon of `anomalous spreading speeds', 
which was noted in a related context
in \citet{MR2322849}\index{Alsmeyer, G.}\index{Iksanov, A.}. Certain martingales\index{martingale} were deployed in  \citet{MR0400438}.  These have 
been a fruitful topic in their own right, and have probably received more attention
since then than the first-birth problem itself (see \citet{MR2471666} for a recent nice contribution on when these martingales are integrable). However, those developments will be ignored here.

\section{The basic model}
The branching random walk (BRW)  
starts with a single particle
located at the
origin. This particle produces 
daughter particles, which are
scattered in $\mbR$, to give the first generation. These first
generation particles produce their own daughter particles 
similarly to give the sec\-ond
generation, and so on.
Formally, each family is described by 
the collection of points in
$\mbR$  giving the positions of the daughters relative to the parent. 
Multiple points are allowed, so that
in a family there may be
several daughter particles born in the same place.
As usual in branching processes,
the $n$th generation particles reproduce independently
of each other. The process is assumed
supercritical\index{branching process!supercritical branching process}, 
so that the expected family size exceeds one (but need not be finite---indeed even the family size itself need not be finite). 
Let $\tP$ and $\tE$ be
the probability and expectation for this process and 
let $Z$ be the generic reproduction process of points in $\mbR$.
 Thus,
$\tE Z$ is the intensity measure\index{intensity measure|(} of $Z$ and $Z(\mbR)$ is the family size,
which will also be written as $N$. The assumption that the process is supercritical becomes
 that $\tE Z(\mbR)=\tE N>1$. 
To avoid burdening the description with qualifications about the survival set, 
let $\tP (N=0)=0$, 
so that the process survives almost surely. 

The model includes several others. One 
is when each daughter receives an independent displacement,
another is when all daughters receive the 
same displacement, with the 
distribution of the displacement being independent of family size in both cases.
These will be called the BRW with \textit{independent} and \textit{common} displacements respectively. Obviously, in both of these any line of descent follows a trajectory of a random walk. (It is possible to consider an intermediate case, where displacements have these properties conditional on family size, but that is not often done.) 
\label{brw}
Since family size and displacements are independent, these two processes can be coupled in a way that shows that results for one will readily yield results for the other.
In a common displacement BRW imagine each particle occupying the (common) position of its family. Then the process becomes an independent displacement BRW, with a random origin given by the displacement of the first family, and its $n$th generation occupies the same positions as the $(n+1)$th generation in the original common displacement BRW. Really this just treats each family as a single particle.

In a different direction, the points of $Z$ can be confined to $(0,\infty)$
and interpreted as the mother's age at the birth of that daughter: the
framework adopted in
\citet{MR0400438}\index{branching process!age-dependent branching process|)}\index{Kingman, J. F. C.}. Then the process is the 
general branching process\index{branching process!general branching process} associated with the
names of Ryan\index{Ryan Jr., T.},
Crump\index{Crump, K. S.}, Mode\index{Mode, C. J.} and Jagers\index{Jagers, P.}. 
Finally, when all daughters receive the same positive displacement with a distribution
independent of family size the process is the \index{Bellman, R. E.!Bellman--Harris process}Bellman--Harris branching process: the framework
adopted in \citet{MR0370721}\index{Hammersley, J. M.|(}.

There  are other `traditions', which consider the BRW but introduce and describe it rather differently and usually with other problems in focus. 
There is a long tradition phrased in terms of `\index{multiplicative cascade}multiplicative cascades' 
(see for example \citet{MR1741808}\index{Liu, Q.} and the references there) and a rather
shorter one phrased in terms of `weighted
branching'\index{branching process!weighted branching process} (see for example \citet{MR2199054} and the references there). 
The model has arisen in one form or another in a variety of areas. The 
most obvious is as a model for a population 
spreading through an homogeneous habitat. It has also arisen
in modelling random fractals\index{fractal}
 \citep{MR1785625} commonly in the language of multiplicative cascades,
in the  theoretical study of algorithms\index{algorithm}
\citep{MR1140708}, 
in a problem in \Index{group theory} \citep{MR2114819}
and as an ersatz
for both lattice-based models of spin glasses\index{spin glass} in physics 
\citep{MR1601733} 
and a \Index{number theory} problem
\citep{MR1143401}.

\section{Spreading out: old results}\index{spreading out|(}

 Let $Z^{(n)}$ be the positions occupied by the $n$th generation and
 $\rBo{n}$
 its rightmost point,
 so that
  \[
\rBo{n}= \sup\{z: z \mbox{~a point of~} Z^{(n)}\}.
 \]
One can equally well consider the leftmost particle, and the earliest studies
did that. Reflection of the whole process around the origin
shows the two are equivalent: all discussion here will be expressed
in terms of the rightmost particle. 
The first result, stated in a moment, concerns $\rBo{n}/n$ converging to a constant,
$\Gamma$, which can  reasonably be interpreted as the speed of spread in the positive direction.

\label{spreading out}
A critical role in the theory is played by the Laplace
transform\index{Laplace, P.-S.!Laplace transform} of the intensity measure\index{intensity measure|)} $\tE Z$: 
let $\kappa(\phi)= \log  \int e^{\phi z} \tE Z(dz)$ for $\phi \geq 0$ and $\kappa(\phi)=\infty$
for $\phi<0$. It is easy to see that when this is finite for  some $\phi>0$ the intensity measures of $Z$ and $Z^{(n)} $
are finite on bounded sets, and decay exponentially in their right tail.
The behaviour of the leftmost particle is governed by the behaviour of the transform for negative values of its argument. The definition of $\kappa$ discards these, which simplifies later formulations by automatically 
keeping attention on the right tail. In order to
give one of the key formulae for $\Gamma$ and for later explanation,
 let $ \Fd{\kappa}$ be the  Fenchel dual\index{Fenchel, M. W.!Fenchel dual}  of  $\kappa$,
 which is the convex\index{convexity} function given by
 \begin{equation}\label{Fenchel dual}
 \Fd{\kappa}(a)=\sup_\theta \{\theta a-\kappa(\theta)\}.
 \end{equation}
This is sufficient notation to give the first result.

\begin{theorem}\label{first theorem}
  When there is a $\phi>0$ such that
 \begin{equation}\label{m good}
\kappa(\phi) < \infty,
 \end{equation}
 there is a constant $\Gamma$  such that
  \begin{equation}\label{limit}
 \frac{\rBo{n}}{n}\rightarrow \Gamma\mbox{~~~a.s}.
 \end{equation}
and $ \Gamma= \sup\{a: \Fd{\kappa}(a)<0\}=\inf\{\kappa(\theta)/\theta: \theta\}$.
\end{theorem}
This result was proved for the common BRW with only negative displacements with convergence in probability in 
\citet[Theorem 2]{MR0370721}\index{Hammersley, J. M.|)}.  It was proved
in \citet[Theorem 5]{MR0400438}\index{Kingman, J. F. C.} for $Z$ concentrated on
$(-\infty,0)$ and with $0<\kappa(\phi)<\infty$ instead of (\ref{m good}).
The result stated above is contained in  \citet[Theorem 4]{MR0420890},
which covers the irreducible
multi\-type\index{branching random walk (BRW)!multitype BRW} case also, of which more later.  The second of the formulae
for $\Gamma$ is certainly well-known but cannot be found in the papers
mentioned---I am not sure where it first occurs. 
It is not hard to establish from the first one using the definition and properties of $\Fd{\kappa}$.

The developments described here draw on features of transform theory, to give
properties of $\kappa$, and of convexity theory, to give properties of $\Fd{\kappa}$
and the speed $\Gamma$. There are many presentations of, and notations for, these, tailored to the particular problem under consideration. In this review, results will simply
be asserted. The first of these provides a context for the next theorem and aids interpretation  of  $\sup\{a: \Fd{\kappa}(a)<0\}$ in the previous one. It is that 
when $\kappa$
is finite somewhere on $(0,\infty)$, $\Fd{\kappa}$ is an increasing, convex\index{convexity} function,
which is continuous from the left, with minimum value $-\kappa(0)=-\log  \tE N$, which is less than zero.

A slight change in focus
derives Theorem \ref{first theorem} from the asymptotics of
 the numbers of particles in suitable half-infinite intervals.
As part of the derivation of this the asymptotics of the expected numbers are obtained. Specifically, it is shown that when (\ref{m good}) holds
\[
n^{-1} \log 
\left(\tE Z^{(n)}[na,\infty)
\right) \rightarrow - \Fd{\kappa}(a)
\]
(except, possibly, at one $a$).
The trivial observation that when the expectation
of integer-valued variables decays geometrically the variables themselves must ultimately be zero implies that $\log 
 Z^{(n)}[na,\infty)$ is ultimately infinite on
$\{a: \Fd{\kappa}(a)>0\}$.  
This motivates introducing a notation for 
sweeping positive values of $\Fd{\kappa}$, and later other functions, 
to infinity and so we let
\begin{equation}\label{sweep}
 \sweep{f}(a)=\left\{\begin{array}{ll}f(a) &\mbox{when~}f(a)\leq 0\\
\infty&\mbox{when~} f(a)> 0
\end{array}\right. 
 \end{equation}
and $\swFd{\kappa}=\sweep{(\Fd{\kappa})}$. The next result can be construed as saying that in crude asymptotic terms this is the only
way actual numbers differ from their expectation.
\begin{theorem}\label{second theorem} When\/ \eqref{m good} holds,
\begin{equation}\label{describe numbers}
\frac{1}{n} \log 
\left(Z^{(n)}[na,\infty)
\right) \rightarrow - \swFd{\kappa}(a) \mbox{~~~a.s.},
\end{equation}
for all $a \neq \Gamma$.
\end{theorem}
\noindent
From this result, which is \citet[Theorem  2]{MR0464415}, 
and the properties of $\swFd{\kappa}$, Theorem \ref{first theorem}
 follows directly.

A closely related continuous-time model arises when the temporal development
is a Markov branching process\index{Markov, A. A.!Markov branching process}
(Bellman--Harris\index{Bellman, R. E.!Bellman--Harris process} with
exponential lifetimes) or even a Yule process\index{Yule, G. U.!Yule process}
(binary splitting too) and movement is Brown\-ian, giving binary branching
Brown\-ian motion\index{Brown, R.!branching Brownian motion}. The process
starts with a single particle at the origin, which then moves with a
Brown\-ian motion with variance parameter $V$. This particle splits in two at
rate $\lambda$, and the two particles continue, independently, in the same way
from the splitting point. (Any discrete skeleton\index{skeleton!discrete skeleton} of
this process is a branching random walk.)

Now, let $\rBo{t}$ be the position of the rightmost particle at time $t$.  Then $
u(x,t)=\tP(\rBo{t} \leq x)$
satisfies  the
\index{Fisher, R. A.!Fisher/KPP equation}(Fisher/Kolmogorov--Petrovski--Piscounov) equation 
\begin{equation}\label{F-KPP}
\frac{\partial u}{\partial t}=V\frac{1}{2}\frac{\partial ^{2}u}{\partial x^{2}
} - \lambda u(1-u),
\end{equation}
which is easy to see informally by conditioning on what happens in $[0, \delta
t]$. The deep studies of \citet{MR0494541,MR705746}\index{Bramson, M. D.} show, among other things,
that (with $V=\lambda=1$)  $\rBo{t}$ converges in distribution when centred on its \Index{median} and that median is (to $O(1)$)
\[
\sqrt{2}t-\frac{1}{\sqrt{2}}\left(\frac{3}{2} \log  t\right) ,
\]
which implies that $\Gamma=\sqrt{2}$ here.
For the skeleton at integer times,
$\kappa(\theta)=\theta^2/2+1$ for $\theta\geq 0$, and using Theorem \ref{first theorem}  on this confirms that $\Gamma=\sqrt{2}$. Furthermore, for later reference, note that $\theta \Gamma -\kappa(\theta)=0$ when $\theta=\sqrt{2}$. 
\label{bbm}

Theorem \ref{first theorem} is for discrete time, counted by generation. 
There are corresponding results for continuous time, 
where the reproduction
is now governed by a random collection of points in time and space ($\mbR^+ \! \times \mbR $).
The first component gives the mother's age at the birth of this daughter and
the second that daughter's position relative to her mother. Then  the
development in time of the process is that of a
\index{branching process!general branching process}general branching process rather than the
Galton--Watson\index{Galton, F.!Galton--Watson process} development that underpins Theorem \ref{first theorem}. This extension is discussed in \citet{MR1384364} and \citet{MR1601689}. 
In it particles may also move during their lifetime and then branching Brown\-ian motion
becomes a (very) special case.  Furthermore, there are also natural versions of Theorems \ref{first theorem} and \ref{second theorem} when particle positions are in $\mbR^d$ rather than $\mbR$---see \citet[\S4.2]{MR1384364} and references there. 
\section{Spreading out:  first refinements}

Obviously rate-of-convergence\index{rate of convergence} questions follow on from (\ref{limit}). 
An aside in \citet[p33]{MR0433619} noted that, 
typically, $\rBo{n}-n \Gamma$ goes to $-\infty$. 
The following result on this is from \citet[Theorem 3]{MR1629030}, and much of
it is contained also in \citet[Lemma 7.2]{MR1618888}\index{Liu, Q.}. 
 When $\tP(Z(\Gamma,\infty)>0)>0$, so displacements greater than $\Gamma$ are possible,  and (\ref{m good}) holds,
there is a finite $\vartheta>0$ with $\vartheta \Gamma-\kappa(\vartheta) = 0$. Thus the condition here, which will recur in later theorems, is not restrictive.
\begin{theorem} \label{theorem to infinity} If there is a finite $\vartheta>0$ 
with $\vartheta \Gamma-\kappa(\vartheta) = 0$, then
\begin{equation}\label{to infinity}
\rBo{n}-n \Gamma \rightarrow -\infty\mbox{~~~a.s.,}
\end{equation}
and the condition is also necessary when $\tP(Z(\Gamma,\infty)>0)>0$.
\end{theorem}

The theorem leaves some loose ends when $\tP(Z(\Gamma,\infty)=0)=1$.
Then  $\rBo{n}-n \Gamma$ is a decreasing sequence, and so it does have a
limit, but whether (\ref{to infinity}) holds or not is really the explosion
(i.e.\ regularity) problem for the general branching process\undex{branching process!general}:  whether,  with a point $z$ from $Z$ corresponding to a
birth time of $\Gamma-z$, there can be an infinite number of births in a
finite time. This is known to be complex---see \citet{MR0359040}\index{Grey, D. R.} for example.
In the simpler cases it is properties of  $Z(\{\Gamma\})$, 
the number of daughters displaced by exactly $\Gamma$, that matters.

If $Z(\{\Gamma\})$ is the family size of a surviving
branching process (so either $\tE Z(\{\Gamma\}) >1$ or $\tP(Z(\{\Gamma\})=1)=1$) 
it is easy to show that $(\rBo{n}-n \Gamma)$
has a finite limit---so (\ref{to infinity}) fails---using
embedded surviving processes resulting from
focusing on daughters displaced by $\Gamma$: 
 see \citet[Proposition II.5.2]{my-thesis} or \citet[Theorem
 1]{MR1133373}\index{Dekking, F. M.|(}\index{Host, B.|(}. 
In a similar vein, with extra conditions, \citet[Theorem
4]{Addarioberryreed}\index{Addario-Berry, L.}\index{Reed, B.|(} show $\tE(\rBo{n}-n \Gamma) $ is  bounded.

Suppose now that (\ref{m good}) holds. 
When  $\tP(Z(a,\infty)=0)=1$, simple properties of
transforms imply that $\theta a-\kappa(\theta) \uparrow - \log  \tE Z(\{a\}) $ 
as $\theta \uparrow \infty$. 
Then, when $\tE Z(\{a\}) <1$
a little \Index{convexity} theory shows that $\Gamma<a$ and that 
there is a finite $\vartheta$ with $\vartheta \Gamma-\kappa(\vartheta) = 0$, 
so that Theorem \ref{theorem to infinity} applies. 
This leaves the case where 
(\ref{m good}) holds, $\tP(Z(\Gamma,\infty)=0)=1$ and  $\tE Z(\{\Gamma\}) =1$ but $\tP( Z(\{\Gamma\})=1)<1$, 
which is sometimes called,
misleadingly in my opinion, the \textit{critical} branching random
walk\undex{branching random walk (BRW)!critical BRW}
because the process of daughters displaced by exactly $\Gamma$ from their
parent forms a critical Galton--Watson process\index{Galton, F.!Galton--Watson
process}. For this case,
\citet[Theorem 1]{MR510529}\index{Bramson, M. D.|(} and \citet[\S9]{MR1133373} show
that (\ref{to infinity}) 
holds under extra conditions including that displacements lie in a lattice, 
and that the convergence is at  rate $\log  \log  n$.  
 \citet[Theorem 2]{MR510529}
also gives conditions under which (\ref{to infinity}) fails.

\section{Spreading out:  recent refinements}

The challenge to derive analogues for the branching random walk
 of the fine results for
branching Brown\-ian motion\index{Brown, R.!branching Brownian motion} 
has been open for a long time. 
Progress was made
in \citet{MR1325045}\index{McDiarmid, C.|(} and, very recently, a nice result has been given in
\citet[Theorem 1.2]{hu-shi}\index{HuY@Hu, Y.}\index{ShiZ@Shi, Z.}, under reasonably mild conditions.
Here is its translation into the  current notation. It shows that the numerical identifications noted in the branching Brown\-ian motion case 
in \S\ref{bbm}  are general.
\begin{theorem}\label{h-s}
Suppose that there is a $\vartheta>0$ with  $\vartheta \Gamma-\kappa(\vartheta) = 0$, 
and that, for some $\epsilon>0$, $ \tE(N^{1+\epsilon})<\infty$,
$\kappa(\vartheta+\epsilon)<\infty$ and $\int e^{-\epsilon z}\tE Z(dz)<\infty$. 
Then
\[
-\frac{3}{2}=\liminf_n \frac{\vartheta(\rBo{n}-n \Gamma)}{ \log  n}
<\limsup_n \frac{\vartheta(\rBo{n}-n \Gamma)}{ \log  n}=-\frac{1}{2}\mbox{~~~a.s.}
\]
and
\[
 \frac{\vartheta(\rBo{n}-n \Gamma)}{ \log  n} \rightarrow -\frac{3}{2}\mbox{~~in probability}.
\]
\end{theorem}
Good progress has also been made on the tightness\index{tight|(} of the distributions of $\rBo{n}$ when centred
suitably.  Here is a recent  result from \citet[Theorem
1.1]{Bramsonzeitouni}\index{Bramson, M. D.|)}\index{Zeitouni, O.}.
\begin{theorem}Suppose the BRW has independent or common displacements according to the random variable $X$. Suppose also that for some $\epsilon>0$, $ \tE(N^{1+\epsilon})<\infty$ and that for some $\psi>0$ and $y_0>0$
\begin{equation}\label{BZ}
\tP(X>x+y)\leq e^{-\psi y}\tP(X>x) ~~~~\forall x>0, y>y_0.
\end{equation}
Then the distributions of $\{\rBo{n}\}$ are tight when centred on their medians.
\end{theorem}
\noindent
It is worth noting that (\ref{BZ}) ensures that (\ref{m good}) holds for all $\phi \in [0,\psi)$. 
There are other results too---in particular, \citet[Theorem
1]{MR1325045}\index{McDiarmid, C.|)} and  \citet[\S3]{MR1133373}\index{Dekking, F. M.|)}\index{Host, B.|)} both give tightness results for the (general) BRW, 
but with $Z$ concentrated on a half-line. Though rather old for this section,  \citet[Theorem 2]{MR1133373} is worth recording here: the authors assume the BRW is concentrated on a lattice, but they do not use that in the proof of this theorem. To state it, let 
$\widetilde{D}$ be the second largest point in $Z$ when 
$N\geq 2$ and the only point otherwise. 
\begin{theorem} If the points of $Z$ are confined to $(-\infty,0]$ and  $\tE \widetilde{D}$ is finite, then $\tE \rBo{n}$ is finite and the distributions of $\{\rBo{n}\}$ are tight when centred on their expectations.
\end{theorem}\noindent
The condition that $\tE \widetilde{D}$ is finite holds when $\int e^{\phi z}\tE Z(dz)$ is finite in a neighbourhood of the origin, which is contained within the conditions in Theorem \ref{h-s}.
In another recent study  \citet[Theorem 3]{Addarioberryreed}\index{Addario-Berry, L.}\index{Reed, B.|)} give the following result, 
which gives tightness and also estimates the centring. 
\begin{theorem}\label{ab-r}
Suppose that there is a $\vartheta>0$ with  $\vartheta \Gamma-\kappa(\vartheta) = 0$, 
and that, for some $\epsilon>0$, 
$\kappa(\vartheta+\epsilon)<\infty$ and $\int e^{-\epsilon z}\tE Z(dz)<\infty$. 
  Suppose also that the BRW has a 
finite maximum family size and  independent displacements.
Then
\[
\tE\rBo{n}=n \Gamma  -\frac{3}{2 \vartheta}\log  n+O(1),
\]
and there are $C>0$ and $\delta>0$ such that
\[
\tP\left(|\rBo{n}-\tE\rBo{n}|>x\right) \leq Ce^{-\delta x} ~~~\forall x.
\]
\end{theorem}\noindent
The conditions in the first sentence here have been stated in a way that keeps
them close to those in Theorem \ref{h-s} rather than specialising them for
independent displacements.   Now, moving from tightness\index{tight|)} to convergence in
distribution---which cannot be expected to hold without a non-lattice
assumption---the following result, which has quite restrictive conditions, is taken from
 \citet[Theorem 1]{MR1765165}\index{Bachmann, M.}. 
\begin{theorem}\label{Bach}
Suppose that the BRW has $ \tE N<\infty$ and independent displacements according to a random variable with density function $f$ where $- \log  f$ is convex.
Then the variables $\rBo{n}$ converge in distribution when centred on medians. 
\end{theorem}\index{convexity}\index{median}
It is not hard to use the coupling mentioned in \S\ref{brw} to see that Theorems \ref{ab-r} and \ref{Bach} imply that these two results also hold for common displacements. 
\section{Deterministic theory}\label{deter}\index{deterministic modelling}
There is another, deterministic, stream of work concerned with
modelling the  spatial \index{spread of population}spread of populations 
in a homogeneous habitat, and
closely linked to the study of reaction-diffusion equations\index{reaction-diffusion equation|(} like (\ref{F-KPP}).  
The main presentation  is 
\citet{MR653463}\index{Weinberger, H. F.|(}, with a formulation that 
has much in common with that adopted in \citet{MR0370721}\index{Hammersley,
J. M.}. Here the description of the framework is pared-down.
This sketch draws heavily on \citet{MR1943224},
specialised to the homogeneous (i.e.\ aperiodic) case and 
one spatial dimension. 
The aim is to say enough to make certain connections with the BRW.

Let $u^{(n)}(x)$
 be the density of the population (or the
\index{gene!frequency}gene frequency, in an alternative interpretation) at time $n$ and position $ x \in \mbR$.  This is a discrete-time theory, 
so there is an
\index{operator!updating operator}updating operator $Q$ satisfying $u^{(n+1)}=Q(u^{(n)})$. More formally, 
let ${\cal F}$ be the non-negative continuous functions on $\mbR$ 
bounded  by $\beta$. Then $Q$ 
maps ${\cal F}$  into itself and
$u^{(n)}=Q^{(n)}(u^{(0)})$, where $u^{(0)}$ 
is the initial density and $Q^{(n)}$ is the $n$th iterate of $Q$. 
The operator is to satisfy the following restrictions. 
The constant functions at $0$ and at $\beta$ are both fixed points of $Q$. 
For any function $u \in {\cal F}$ that is not zero everywhere, $Q^{(n)}(u) \rightarrow \beta$, and  $Q(\alpha) \geq \alpha$ for non-zero constant functions in ${\cal F}$.  (Of course, 
without the spatial component, this is all reminiscent of the basic properties of the generating function of the family-size.) 
The operator $Q$ is
\index{order preserving@order-preserving}order-preserving,
in that if $u \leq v$ then $Q(u)\leq Q(v)$, so increasing the population anywhere
never has deleterious effects\index{deleterious effect} in the future; it is also \Index{translation-invariant},
because the habitat is homogeneous, and suitably continuous. Finally, every sequence
$u_m \in {\cal F}$ 
contains a subsequence $u_{m(i)}$ such that $Q(u_{m(i)})$ converges uniformly on compacts.
Such a $Q$ can be 
obtained by taking the development of a
\index{reaction-diffusion equation|)}reaction-diffusion equation for a time $\tau$. Then
$Q^{(n)}$ gives the development to time $n \tau$, and the results for this discrete
formulation transfer to such equations.

Specialising \citet[Theorem 2.1]{MR1943224}, there is a spreading speed  $\Gamma$ in the following sense.
If $u^{(0)}(x)=0$ for $x\geq L$ and $u^{(0)}(x)\geq \delta>0$ for all $x\leq K$,  then for any $\epsilon >0$ 
\begin{equation}\label{spreading}
\sup_{x\geq n(\Gamma+\epsilon)} |u^{(n)}(x)| \rightarrow 0 \mbox{~~and~~}
\sup_{x\leq n(\Gamma-\epsilon)} |u^{(n)}(x)-\beta|\rightarrow 0.
\end{equation}
In some cases the spreading speed can be computed through
linearisation\break---see \citet[Corollary 2.1]{MR1943224}\index{Weinberger, H. F.|)} and \citet[Corollary to
Theorem 3.5]{Lui1989269}\index{Lui, R.}---in that the speed is the same as that obtained by 
replacing $Q$ by a truncation of its linearisation at the zero function. So
 $Q(u)=Mu$ for small $u$ and $Q(u)$ is replaced by $\min\{\omega, Mu\}$, where $\omega$ is a constant, positive function with $M \omega > \omega$.
The linear functional $Mu(y)$ must be represented as an integral with respect to some measure, 
and so, using the translation invariance of $M$, there is a measure $\mu$ such that
\begin{equation}\label{mu}
Mu(y)=\int u(y-z)\,\mu(dz).
\end{equation}
Let $\tilde{\kappa}(\theta)=\log  \int e^{\theta z} \mu(dz)$. Then the results show that 
the speed $\Gamma$ in  (\ref{spreading}) is given by
\begin{equation}\label{W-Gamma}
\Gamma = \inf_{\theta>0} \frac{\tilde{\kappa}(\theta) }{\theta}.
\end{equation}
Formally, this is one of the formulae for the speed in Theorem \ref{first theorem}.  In fact, the two frameworks can be linked, as indicated next.

In the BRW, suppose the generic reproduction process $Z$ has points $\{z_i\}$. 
Define $Q$
by
\[
Q\left(u(x)\right)=1-\tE \left[1-\prod_i u (x-z_i)\right].
\]
This has the general form described above with $\beta=1$.
On taking $u^{(0)}(x)=\tP (\rBo{0}>x)$ (i.e.\ Heaviside\index{Heaviside, O.} initial data) it is easily established by induction that 
$
u^{(n)}(x)=\tP (\rBo{n}>x)$.
This is in essence the same as the observation that the 
distribution of the rightmost particle in branching Brown\-ian
motion\index{Brown, R.!branching Brownian motion} satisfies the differential equation (\ref{F-KPP}). 
The idea is explored in the  spatial spread of
the `deterministic simple epidemic'\index{epidemic!deterministic simple epidemic}
in \citet{Mollison1993147}\index{Mollison, D.}\index{Daniels, H. E.}, 
a continuous-time model which, like branching Brown\-ian motion,
has BRW as its discrete skeleton\index{skeleton!discrete skeleton}. 
Now Theorem \ref{first theorem}
implies that (\ref{spreading}) holds, and that, for $Q$ obtained in this way, the speed is indeed given by the (truncated) linear approximation. The other theorems about $\rBo{n}$ also translate into results about such $Q$. For example,  Theorem \ref{Bach} gives conditions for $u^{(n)}$ when centred suitably to converge to a fixed (\Index{travelling wave}) profile.

\section{The multi\-type case}
\label{multitype}\index{branching random walk (BRW)!multitype BRW|(}

 Particles now
have types\index{type} drawn from a finite set,
$\cS$,  and
their reproduction is defined by random points in
$\cS \times \mbR$. The distribution of these points
depends on the parent's type. 
The first component gives the daughter's type
and the second component gives the
daughter's position, relative to the parent's.
As previously, $Z$ is the generic reproduction process,
but now let $Z_\sigma$
be the points (in $\mbR$) corresponding to
 those of type $\sigma$;
$Z^{(n)}$ and $Z^{n}_{\sigma}$ are defined similarly.
Let $\tP_{\!\nu}$ and $\tE_\nu$ be
the probability and expectation associated
with reproduction from an initial ancestor with type $\nu \in \cS$.
Let
 $\rB{\sigma}{n}$
 be the rightmost particle of type $\sigma$ in the $n$th generation,
 and let $\rBo{n}$ be the rightmost of these, which is consistent with 
the one-type notation.

The type space\index{type!type space} can be classified, using the
relationship `can have a descendant of this type', or, equivalently, using the
non-negative expected family-size matrix\undex{branching random walk (BRW)!expected family-size matrix}, $\tE_\nu  Z_\sigma (\mbR)$.
Two types are in the same class when each can have a descendant 
of the other type in some generation. 
When there is a single class 
the family-size matrix is
\textit{irreducible}\index{branching random walk (BRW)!irreducible BRW} and the process is similarly described.
When the expected family-size matrix is \textit{aperiodic} (i.e.\ primitive)
the process is also called
aperiodic\undex{branching random walk (BRW)!aperiodic BRW|(}, and it is
supercritical\undex{branching random walk (BRW)!supercritical BRW|(} when this matrix has Perron--Frobenius\index{Perron, O.!Perron--Frobenius eigenvalue} (i.e.\ non-negative and of maximum modulus) eigenvalue
  greater than one. 
Again, to
avoid qualifications about the survival set, assume
extinction\index{extinction} is impossible from the starting type used.

For $\theta \geq 0$,
let $\exp(\kappa(\theta))$ be the Perron--Frobenius eigenvalue of the matrix of transforms $\int e^{\theta z} \tE_\nu Z_\sigma(dz)$, and let $\kappa(\theta)=\infty$ for $\theta<0$.
If there is just one type, this definition agrees with that of $\kappa$ at the start of \S\ref{spreading out}.
The following result, which is \citet[Theorem 4]{MR0420890}, has been mentioned already. 
\begin{theorem}\label{multitype first theorem} Theorem\/ $\ref{first theorem}$
holds for any initial type in a supercritical irreducible BRW. 
\end{theorem}
\noindent
The simplest multi\-type version of Theorem \ref{second theorem} is the following, which is proved in \cite{JDB-anom}. When $\sigma=\nu$ it is a special case of results indicated in \citet[\S4.1]{MR1601689}.

\begin{theorem} \label{supercrit}  For a supercritical aperiodic BRW for which\/ \eqref{m good} holds,
  \begin{equation}\label{exp growth ub 2}
\frac{1}{n}\log  \left(Z^{(n)}_\sigma[na, \infty) \right) \rightarrow  - \swFd{\kappa}(a)   \mbox{~~~a.s.-}\tP_{\!\nu}
 \end{equation}for $a \neq \sup\{a:\Fd{\kappa}(a)<0\}=\Gamma$, and
 \[
 \frac{\rB{\sigma}{n}}{n} \rightarrow \Gamma 
\mbox{~~~a.s.-}\tP_{\!\nu}.
 \]
\end{theorem}\undex{branching random walk (BRW)!aperiodic BRW|)}\undex{branching random walk (BRW)!supercritical BRW|)}
\noindent
Again there is a deterministic\index{deterministic modelling} theory, following the pattern described
in \S\ref{deter} and  discussed in \citet{Lui1989269, Lui1989297}\index{Lui, R.}, which can be related to Theorem \ref{multitype first theorem}. 
Recent developments in that area raise some interesting questions that are the subject of the next two sections.

\section{Anomalous spreading}\index{anomalous spreading|(}

In the multi\-type version of the deterministic context of \S\ref{deter},
 recent papers \citep{MR1930974, MR2322849,MR1930975,Li200582} 
have considered what happens 
when the type space\index{type!type space} is reducible. 
 Rather than set out the framework in 
its generality,  the simplest possible case,
the reducible\index{branching random walk (BRW)!reducible BRW|(}\index{branching random walk (BRW)!two-type BRW|(} two-type\index{type!two-type} case, will
be considered here, for the principal issue can be illustrated
through it. The two types will be $\nu$ and $\eta$.
Now, 
the vector-valued non-negative function $u^{(n)}$ gives the
\Index{population density} of two species---the two types, $\nu$ and $\eta$---at $x \in \mbR$ at time $n$, 
and $Q$ models growth, interaction and migration\undex{branching random walk (BRW)!growth, interaction and migration}, as the populations develop in discrete time. The programme is the same as that indicated
in \S\ref{deter}, that is to investigate the existence of spreading speeds
\index{spreading speed|(}and when these speeds can be obtained from the truncated linear approximation\undex{branching random walk (BRW)!truncated linear approximation}. 

In this case the approximating
\index{operator!linear operator}linear operator, generalising that given in (\ref{mu}), is
\begin{align*}
(Mu(y))_\eta&=\int u_\eta (y-z) \,\mu_{\eta \eta}(dz),\\
(Mu(y))_\nu&=\int u_\nu (y-z)\,\mu_{\nu \nu} (dz)
+\int u_\eta (y-z)\,\mu_{\nu  \eta } (dz).
\end{align*}
Simplifying even further, assume there is no spatial spread 
associated with the `interaction' term here, so that $\int u_\eta (y-z)\,\mu_{\nu\eta} (dz)= c u_\eta (y)$ for some $c>0$.
The absence of $\mu_{ \eta \nu}$ in the first of these
makes the linear approximation reducible. The first equation 
is really just for the type $\eta$ and so  will have the speed that corresponds to
$\mu_{\eta \eta}$, given through its transform by (\ref{W-Gamma}), and written $\Gamma_\eta$. In the second, on ignoring the interaction term, 
it is plausible that the speed must be at least that of type $\nu$ alone, which corresponds to
$\mu_{\nu \nu}$ and is written $\Gamma_\nu$. However, it can also have
 the speed of $u_\eta$ from the `interaction' term.
It is claimed in \citet[Lemma 2.3]{MR1930974}\index{Weinberger,
H. F.|(}\index{Lewis, M. A.|(}\index{LiB@Li, B.|(} that 
when $Q$ is replaced by the approximating operator 
$\min\{\omega, Mu\}$ this does behave as just outlined, with the corresponding formulae for the speeds: thus that of $\eta$ is $\Gamma_\eta$ and that for 
$\nu$ is $\max \{\Gamma_\eta,\Gamma_\nu\}$.
However, in  \citet{MR2322849} a flaw in the argument is noted, and an example is given where the speed of $\nu$ in the truncated linear approximation can be faster than this, 
the anomalous spreading speed of their title, though the actual speed is not identified.   The relevance of the phenomenon to a biological example 
is explored in
  \citet[\S5]{MR2322849}.

As in \S\ref{deter}, the BRW provides some particular examples of $Q$ that fall
within the general scope of the deterministic theory.
Specifically, suppose the generic reproduction process $Z$ has points
$\{\sigma_i,z_i\}\in \cS \times \mbR$. Now let $Q$, which operates on vector
functions indexed by the type space\index{type!type space} $\cS$, be defined by
\[
Q\left(u(x)\right)_{\nu}=1-\tE_\nu \left[1-\prod_i u_{\sigma_i}(x-z_i)\right].
\]
Then, just as in the one-type case, when $u^{(0)}_{\nu}(x)=\tP_{\!\nu} (\rBo{0}>x)$ induction establishes that $u^{(n)}_{\nu}(x)=\tP_{\!\nu} \left(\rBo{n}>x\right)$. It is perhaps 
worth noting that in the BRW the index $\nu$ is the starting type, whereas it is the `current' type in \citet{MR2322849}.
However, this makes no formal difference.

Thus, the anomalous spreading phenomenon should be manifest in the BRW,
and, given the more restrictive framework, it should be possible to pin down 
the actual speed there, and hence for the corresponding $Q$ with
Heaviside\index{Heaviside, O.} initial data. 
This is indeed possible. 
Here the discussion stays with the simplifications already used in looking at the deterministic results.

Consider a two-type BRW in which each type $\nu$ always produces at least one daughter of type $\nu$, on average produces more than one, and can produces daughters of
type $\eta$---but type $\eta$ never produce daughters of  type $\nu$. Also for $\theta \geq 0$ let
\[
\kappa_{\nu}(\theta)=\log  \int e^{\theta z} \tE_\nu Z_\nu(dz) \mbox{~~and~~}\kappa_{\eta}(\theta)=\log  \int e^{\theta z} \tE_\eta Z_\eta(dz)
\]
and let these be infinite for $\theta < 0$.
Thus Theorem  \ref{second theorem} applies to type $\nu$ considered alone to show that
\[
\frac{1}{n}\log  \left(Z^{(n)}_\nu[na, \infty) \right) \rightarrow  - \swFd{\kappa}_\nu(a)   \mbox{~~~a.s.-}\tP_{\!\nu}.
\]
It turns out that this estimation of numbers is critical in establishing the speed 
for type $\eta$. It is possible
 for the growth in numbers of type $\nu$, through the numbers
 of type $\eta$ they produce, to increase the speed of
type $\eta$ from that of a population without type $\nu$.
This is most obvious if type $\eta$ is subcritical, so that any
line of descent from a type $\eta$ is finite, 
for the only way they can then spread is 
through the `forcing' from type $\nu$.
However, if in addition the dispersal distribution at reproduction 
for $\eta$ has a much heavier tail than that for $\nu$
it is now possible for type $\eta$  to spread faster than
type $\nu$.

For any two functions $f$ and $g$, let $\cv{f}{g}$
be the greatest (\Index{lower semi-continuous}) convex\index{convexity|(} function
beneath both of them. The following result is a very special case 
of those proved in \citet{JDB-anom}.  
The formula given in the next result for the  speed $\Gamma^{\dagger}$ is the same as that given in \citet[Proposition 4.1]{MR2322849} as the upper bound on the speed of the truncated  linear approximation.

\begin{theorem} \label{prelim main theorem} Suppose 
 that $\max\{\kappa_{\nu}(\phi_\nu), \kappa_{\eta}(\phi_\eta)\}$ is finite for some $\phi_\eta \geq \phi_\nu >0 $ and that
\begin{equation}\label{off-diag}
\int e^{\theta z} \tE_\nu Z_\eta(dz)<\infty ~~~\forall \theta \geq 0.
\end{equation}
    Let $r=\sweep{\cv{\swFd{\kappa}_\nu}{\Fd{\kappa}_\eta}}$.
Then
\begin{equation}\label{key result}
\frac{1}{n} \log
\left(Z^{(n)}_{\eta}[na,\infty)
\right) \rightarrow -r(a) \mbox{~~~a.s.-}\tP_{\!\nu},
\end{equation}
for  $a \neq  \sup\{a: r(a)<0\}=\Gamma^{\dagger}$, and
\begin{equation}\label{seq speed}
 \frac{\rB{\eta}{n}}{n} \rightarrow \Gamma^{\dagger}. \mbox{~~~a.s.-}\tP_{\!\nu}.
\end{equation}
Furthermore,
\begin{equation}\label{two classes}
 \Gamma^{\dagger}
= \inf_{0<\varphi \leq \theta}
\max\left\{\frac{\kappa_\nu(\varphi)}
{\varphi},\frac{\kappa_\eta(\theta)}
{\theta}
\right\}.
\end{equation}
\end{theorem}

\begin{figure}[!ht]
\includegraphics[scale=0.3]{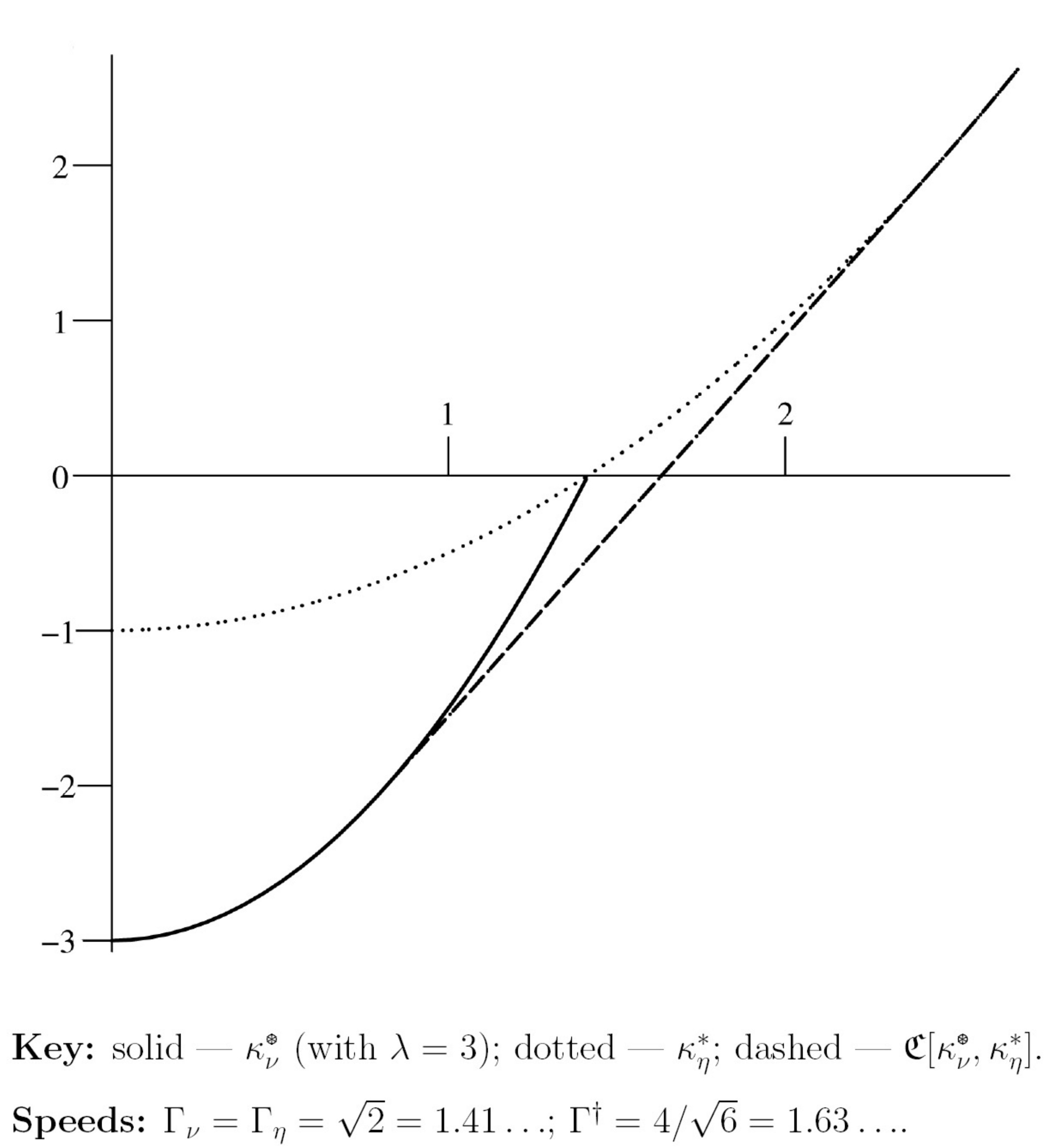}
\caption{Illustration of how anomalous speed arises. }
\label{ff1}
\end{figure}

From this result it is possible to see how  $ \Gamma^{\dagger}$
can be anomalous.
Suppose that
$r( \Gamma^{\dagger})=0$, so that 
$ \Gamma^{\dagger}$ is the speed, and that
$r$ 
is strictly below both $\swFd{\kappa}_\nu$ and $\Fd{\kappa}_\eta$
at  $ \Gamma^{\dagger}$.  This will occur when the minimum of the 
two convex functions  $\Fd{\kappa}_\nu$ and $\Fd{\kappa}_\eta$ is not convex at $\Gamma^{\dagger}$, and then the largest  convex\index{convexity|)} function below both will be linear there. In these circumstances, $\Fd{\kappa}_\nu(\Gamma^{\dagger})>0$, which implies that $\swFd{\kappa}_\nu(\Gamma^{\dagger})=\infty$,
and $\Fd{\kappa}_\eta(\Gamma^{\dagger})>0$. Thus $\Gamma^{\dagger}$
will be strictly greater than both $\Gamma_{\nu}$ and
$\Gamma_{\eta}$, giving a `super-speed'---Figure \ref{ff1} illustrates a case that will soon be described fully 
where $\Gamma_{\nu}$ and $\Gamma_{\eta}$ are equal and $\Gamma^{\dagger }$ exceeds them.
 Otherwise, that is when $\Gamma^{\dagger}$ is not in a linear portion of $r$, $\Gamma^{\dagger}$ is just the maximum of
$\Gamma_{\nu}$ and  $\Gamma_{\eta}$. 

The example in \citet{MR2322849}\index{Weinberger, H. F.|)}\index{Lewis, M. A.|)}\index{LiB@Li, B.|)} that illustrated anomalous speed 
was derived from coupled
\index{reaction-diffusion equation}reaction-diffusion equations. When there is a branching interpretation, 
which it must be said will be the  exception not the rule,
the actual speed can be identified through Theorem \ref{prelim main theorem} 
and its generalisations. 
This will now be illustrated with an example.
Suppose type $\eta$ particles form a binary branching Brown\-ian
motion\index{Brown, R.!branching Brownian motion}, 
with variance parameter and splitting rate both one. 
Suppose type $\nu$ particles form a branching Brown\-ian motion, 
but with variance parameter $V$, splitting rate $\lambda$ and,  
on splitting, type $\nu$ particles produce a (random) 
family of particles of both types. There are $1+N_\nu$ of type  $\nu$ and $N_\eta$ of type $\eta$,  so 
that the family always contains at least one daughter of type $\nu$; the corresponding bivariate
\index{probability generating function (pgf)}probability generating function is $\tE a^{1+N_\nu}b^{N_\eta}=af(a,b)$.
Let $v(x,t)=\tP_\eta(\rBo{t} \leq x)$ and 
$w(x,t)=\tP_\nu(\rBo{t} \leq x)$. These satisfy
\begin{align*}
\frac{\partial v}{\partial t}&=\frac{1}{2}\frac{\partial ^{2}v}{\partial x^{2}
}-v(1-v),\\
\frac{\partial w}{\partial t}&=V\frac{1}{2}\frac{\partial ^{2}w}{\partial x^{2}
}-\lambda w(1-f(w,v)).
\end{align*} 
Here, when  the initial ancestor is of type $\nu$ and at the origin the initial data are  $w(x,0)=1$ for $x \geq 0$ and 0 otherwise and $v(x,0) \equiv 0$. 
Note that, by a simple change of variable, these can be rewritten as 
equations in $\tP_\eta(\rBo{t} > x)$ and 
$\tP_\nu(\rBo{t} > x)$ where the differential parts are unchanged, but the other terms look rather different. 

Now suppose that $af(a,b)=a^2 (1-p+p b)$, so that a type $\nu$
particle always
splits into two type $\nu$ and with probability $p$ also produces one type $\eta$. 
Looking at the discrete skeleton\index{skeleton!discrete skeleton} at integer times, 
$
\kappa_\nu(\theta)=V\theta^2/2+\lambda
$ for $\theta \geq 0$,
giving
\[
\Fd{\kappa}_\nu(a)=\left\{\begin{array}{ll} -\lambda &a < 0
\\
\displaystyle -\lambda+ \frac{1}{2}
\frac{a^2}{V}& a \geq 0
\end{array}\right. 
\]
and speed $(2V\lambda)^{1/2}$, obtained by solving $\Fd{\kappa}_\nu(a)=0$. 
The formulae for 
$\Fd{\kappa}_\eta$ are just the special case with $V=\lambda=1$.  
Now, for convenience, take $V=\lambda^{-1}$, 
so that both types, considered alone, have the same speed.
Then, sweeping positive values to infinity,
\[
\swFd{\kappa}_\nu(a)=
\left\{\begin{array}{ll}-\lambda &a< 0, \\
\displaystyle -\lambda\left(1- 
\frac{a^2}{2}\right)& a \in [0, 2^{1/2}],\\
\infty& a > 2^{1/2}.\end{array}\right. 
\]
Now $\cv{\swFd{\kappa}_\nu}{\Fd{\kappa}_\eta}$ is the largest
convex\index{convexity} function below this and $\Fd{\kappa}_\eta$. When $\lambda=3$ these three functions 
are drawn in Figure \ref{ff1}. 

\noindent The point where 
each of them meets the horizontal axis gives 
the value of speed  for that function. 
Thus, $\Gamma^{\dagger}$
exceeds the other two, which are both $\sqrt{2}$.
Here $\Gamma^{\dagger}=4/\sqrt{6}$. In general, for $\lambda>1$, it is 
$(1+\lambda)/\sqrt{2 \lambda}$,  
which can be made arbitrarily large by increasing $\lambda$ sufficiently.
\section{Discussion of anomalous spreading}

The critical function in Theorem \ref{prelim main theorem} is $r=
\sweep{\cv{\swFd{\kappa}_\nu}{\Fd{\kappa}_\eta}}$. Here is how it arises.
The function $\swFd{\kappa}_\nu$ describes the growth in numbers and spread
of the type $\nu$. Conditional on these, $\cv{\swFd{\kappa}_\nu}{\Fd{\kappa}_\eta}$ describes the growth and spread in expectation of those
of type $\eta$.
To see why this might be so, take a $b$ with $\swFd{\kappa}_\nu(b)<0$ so that (\ref{describe numbers}) describes the exponential growth of $Z^{(m)}_\nu[mb, \infty)$: there are roughly $\exp( -m\swFd{\kappa}_\nu(b))$ such particles in generation $m$. Suppose now, for simplicity, 
that each of these produces a single particle of type $\eta$ at the parent's position. As noted just before Theorem \ref{second theorem}, the expected numbers 
of type $\eta$ particles in  generation $r$ and in $[rc,\infty)$ descended from a single type $\eta$ at the origin is roughly $\exp(-r \Fd{\kappa}_\eta(c))$. Take $\lambda \in (0,1)$ with $m=\lambda n$ and $r=(1-\lambda)n$.  
Then, conditional on  the development of the first $m$ generations, the expectation of the numbers
of type $\eta$ in generation $n$ and to the right of $mb+rc=n (\lambda b+(1-\lambda) c)$ will be (roughly) at least  $\exp(-n(\lambda \swFd{\kappa}_\nu(b)+(1-\lambda) \Fd{\kappa}_\eta(c))) $. As $b$, $c$ and $\lambda$ vary with $\lambda b+(1-\lambda) c=a$, the least value for 
$\lambda \swFd{\kappa}_\nu(b)+(1-\lambda) \Fd{\kappa}_\eta(c)$
is given by $\cv{\swFd{\kappa}_\nu}{\Fd{\kappa}_\eta}(a)$. 
There is some more work to do to show that this lower bound on the conditional expected numbers is also an upper bound---it is here that (\ref{off-diag}) comes into play.
Finally, 
as indicated just before Theorem \ref{second theorem}, 
this corresponds to actual numbers only when negative, so the positive values
of this convex minorant\index{convexity!convex minorant} are swept to infinity.

When the speed is anomalous, this indicative description of how $r=
\sweep{\cv{\swFd{\kappa}_\nu}{\Fd{\kappa}_\eta}}$ arises makes plausible the following description of lines of descent  with speed near $\Gamma^{\dagger}$. They will arise as a `dog-leg', with the first portion of the trajectory, which is a fixed proportion of the whole, being a line of descent of type $\eta$ with a speed less than $\Gamma_{\eta}$. The remainder is a line of descent of type $\nu$, with a speed faster than $\Gamma_{\mu}$ (and also than $\Gamma^{\dagger}$).

Without the truncation, the
\index{operator!linear operator}linear operator
approximating (near $u\equiv1$) a $Q$ associated with a BRW  
describes the development of its expected numbers, and so it is 
tempting to define the speed using this, by look\-ing at when expected numbers start to decay.
In the irreducible case\index{branching random walk (BRW)!irreducible BRW}, 
Theorem \ref{supercrit} has an analogue for expected numbers, that
\[
\frac{1}{n}\log  \left(\tE_{\nu} Z^{(n)}_\sigma[na, \infty) \right) \rightarrow  - \Fd{\kappa}(a),  
\]
and so here the speed can indeed be found by looking at when expected numbers start to decay. In contrast, in the set up in Theorem \ref{prelim main theorem}
\[
\frac{1}{n} \log
\left(\tE_{\nu} Z^{(n)}_{\eta}[na,\infty)
\right) \rightarrow -\cv{\Fd{\kappa}_\nu}{\Fd{\kappa}_\eta}(a), 
\] 
and the limit here can be lower than $\cv{\swFd{\kappa}_\nu}{\Fd{\kappa}_\eta}$---the distinction between the functions is 
whether or not positive values are 
swept to infinity in the first argument. Hence
the speed computed by simply asking when expectations start to decay can be too large. In 
Figure \ref{ff1}, $\cv{\swFd{\kappa}_\nu}{\Fd{\kappa}_\eta}$ 
is the same as $\cv{\Fd{\kappa}_\nu}{\Fd{\kappa}_\eta}$,
but it is easy to see, reversing the roles of ${\Fd{\kappa}_\nu}$ and ${\Fd{\kappa}_\eta}$, that $\cv{\swFd{\kappa}_\eta}{\Fd{\kappa}_\nu}$ is the same as $\Fd{\kappa}_\nu$.
Thus if $\eta$ could produce  $\nu$, rather than the other way round, 
expectations would still give the speed $\Gamma^{\dagger}$
but the true speed would be $\Gamma_\nu
(=\Gamma_\eta)$.\index{branching random walk (BRW)!reducible BRW|)}\index{branching random walk (BRW)!two-type BRW|)}\index{spreading speed|)}

The general case, with many classes, 
introduces a number of additional challenges 
(mathematical as well as notational). It is discussed in \citet{JDB-anom}.  The matrix of transforms now has irreducible blocks on
its diagonal, corresponding to the classes, and their Perron--Frobenius
eigenvalues\index{Perron, O.!Perron--Frobenius eigenvalue}
supply the $\kappa$ for each class, as would be anticipated from \S\ref{multitype}.
Here a flavour of some of the other complications.
The rather strong condition (\ref{off-diag}) means that the spatial distribution 
of type $\eta$ daughters to a type $\nu$ mother 
 is irrelevant to the form of the result. If convergence is assumed only for some $\theta>0$ rather than all this need not remain true. One part of the challenge is to describe when these `off-diagonal' terms remain irrelevant; another is to say what happens when they are not. 
If there are  various routes through the classes
from the initial type to the one of interest these possibilities must be combined: 
in these circumstances, the function $r$ in (\ref{key result})
need not be convex\index{convexity} (though it will be increasing). 
It turns out that
the formula for $\Gamma^{\dagger}$, 
which seems as if it might be particular to the case of two classes,
extends fully---not only 
in the sense that there is a version that involves more classes, 
but also in the sense that 
the speed can usually be obtained as the maximum of that obtained using
(\ref{two classes})  for all pairs of classes where the first can have
descendants in the second (though the line of descent may have to go through
other classes on the way)\index{branching random walk (BRW)|)}\index{spreading out|)}\index{branching random walk (BRW)!multitype BRW|)}\index{anomalous spreading|)}.

\bibliographystyle{cambridgeauthordate}
\bibliography{biblio-kingman}

\end{document}